\documentclass{elsarticle}

\usepackage{hyperref}

\usepackage{color, colortbl}
\definecolor{Gray}{gray}{0.9}

\usepackage{mathrsfs}
\usepackage[binary-units]{siunitx}

\usepackage{stfloats}

\usepackage{pstool}

\usepackage{graphicx,psfrag}

\usepackage{subcaption}

\usepackage{cuted}

\usepackage{mathrsfs}
\usepackage{amsthm}

\newtheorem{theorem}{Theorem}

\newtheorem{lemma}{Lemma}
\newtheorem{problem}{Problem}
\newtheorem{remark}{Remark}

\usepackage{pdfpages}

\usepackage{hyperref} 

\usepackage{algorithm,algorithmic}

\usepackage{dsfont}

\usepackage{amsfonts}
\usepackage{caption}
\usepackage{graphicx}

\usepackage[utf8]{inputenc}
\usepackage{amssymb}

\usepackage{nomencl} 
\makenomenclature
\usepackage{etoolbox}
\renewcommand{\nomgroup}[1]{%
  \item[\bfseries
  \ifstrequal{#1}{A}{Given Parameters}{%
  \ifstrequal{#1}{B}{Decision variables}{%
  \ifstrequal{#1}{C}{Other Symbols}{}}}%
]}

\usepackage{amsmath}
\usepackage{algorithmic}

\usepackage{array}

\usepackage{float}

\usepackage{mathtools}

\usepackage{stfloats}

\newcommand{\pushright}[1]{\ifmeasuring@#1\else\omit\hfill$\displaystyle#1$\fi\ignorespaces}
\newcommand{\pushleft}[1]{\ifmeasuring@#1\else\omit$\displaystyle#1$\hfill\fi\ignorespaces}
\makeatother
\newif\ifmargincomments 
\margincommentstrue

\ifmargincomments

\else

\fi

\ifmargincomments

\else

\fi

\usepackage[a4paper, total={6in, 8in}]{geometry}

\journal{TRC}

\biboptions{authoryear}

\begin{document}

\begin{frontmatter}

\title{Joint Routing and Charging Problem of Electric Vehicles with Incentive-aware Customers Considering Spatio-temporal Charging Prices}



\author[HIT,SUSTech]{Canqi Yao}
\author[SUSTech]{Shibo Chen}
\author[Tue]{Mauro Salazar}
\author[SUSTech]{Zaiyue Yang\corref{mycorrespondingauthor}} \cortext[mycorrespondingauthor]{Corresponding author}
\ead{yangzy3@sustech.edu.cn}


\address[HIT]{School of Mechatronics Engineering, Harbin Institute of Technology, Harbin, 150000, China}
\address[SUSTech]{ Shenzhen Key Laboratory of Biomimetic Robotics and Intelligent Systems, Department of Mechanical and Energy Engineering, and the Guangdong Provincial Key Laboratory of Human-Augmentation and Rehabilitation Robotics in Universities, Southern University of Science and Technology, Shenzhen 518055, China}
\address[Tue]{Control Systems Technology group, Department of Mechanical Engineering, Eindhoven University of Technology, Eindhoven, MB 5600, The Netherlands }


\begin{abstract}

This paper investigates the scheduling problem of a fleet of electric vehicles, providing mobility as a service to a set of time-specified customers, where the operator needs to solve the routing and charging problem jointly for each EV. Hereby we consider incentive-aware customers and propose that the operator offers monetary incentives to customers in exchange for time flexibility. In this way, the fleet operator can achieve a routing and charging schedule with lower costs, whilst the customers receive monetary compensation for their flexibility.
Specifically, we first propose a bi-level optimization model whereby the fleet operator optimizes the routing and charging schedule accounting for the spatio-temporal varying charging price, jointly with a monetary incentive to reimburse the delivery time flexibility experienced by the customers. Concurrently the customers choose their own time flexibility by minimizing their own cost. Second, we cope with the computational burden coming from this nonlinear bi-level optimization model with an accurate reformulation approach consisting of the KKT optimality conditions, a Big-M-based linearization method, and the zero duality gap of convex optimization problems. This way, we convert the proposed problem into a single-level optimization problem, which can be solved by a strengthened generalized Benders decomposition method holding a faster convergence rate than the generalized Benders decomposition method.
To evaluate the effectiveness of the proposed mathematical model, we carry out numerous simulation experiments by using the VRP-REP data of Belgium. The numerical results showcase that the proposed mathematical model can reduce the delivery fees for the customers together with the cost of operation incurred by the fleet operator.
\end{abstract}


\end{frontmatter}

\newpage 


\newpage

\section{Introduction}

\begin{figure}[t]
\centering
\includegraphics[width=0.6\linewidth]{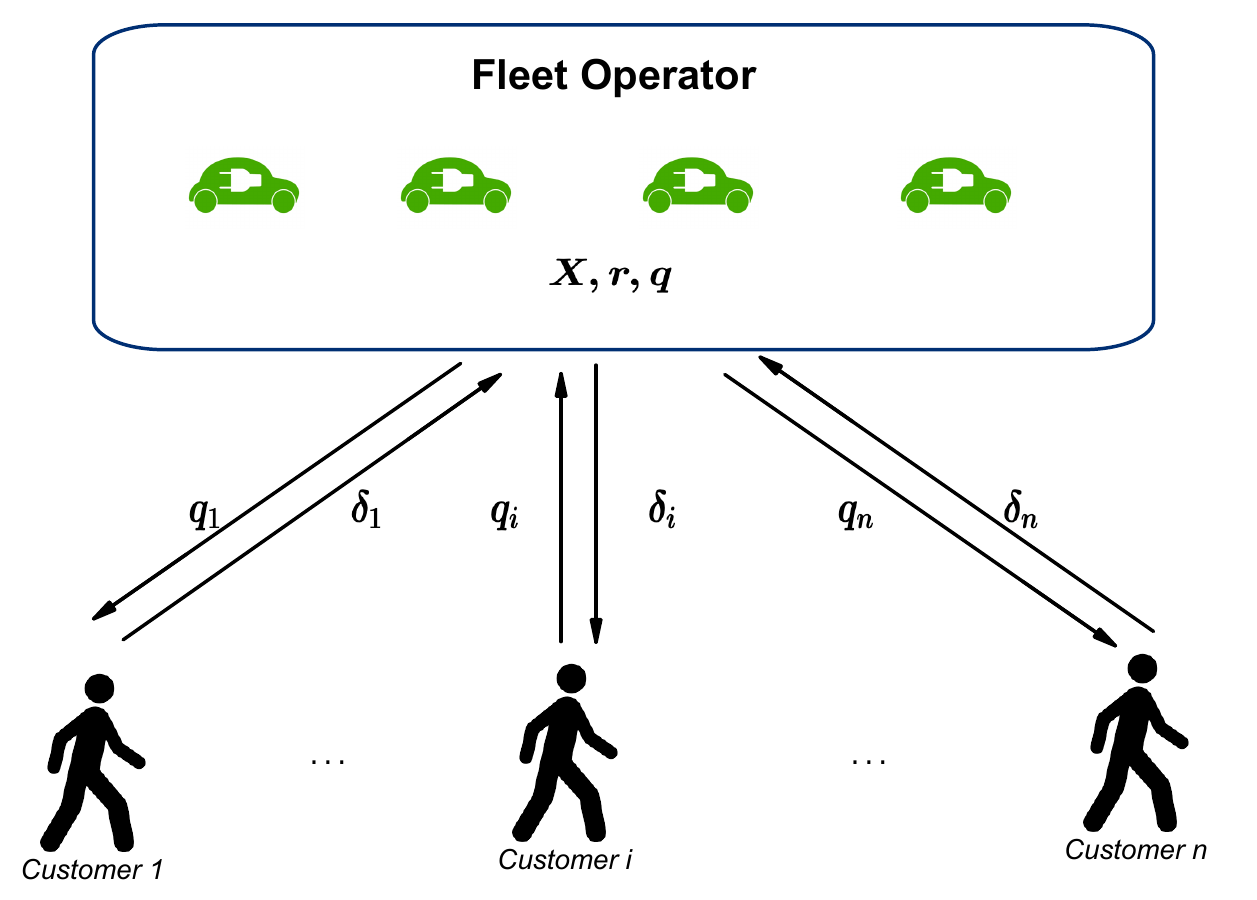}
\caption{The bi-level model of the fleet operator and customers: The fleet operator optimizes the routes of EVs $X$ and the amount of charged energy $r$ jointly with the monetary incentive $q$ offered to the customers to compensate for the time flexibility $\delta$; each customer $i$ chooses the admissible time flexibility $\delta_i$ as a function of the incentive rate $q_i$ provided by the fleet operator.}
\label{bilevel}
\end{figure}

To face a large increase in freight transportation demand while transiting to a carbon neutral society~\cite{carbonneutral}, fleets of electric vehicles (EVs) have been deployed to provide mobility as a service to a set of customers. The fleet operator (refer to Fig.~\ref{bilevel}) aims to make timely decisions on the routing schedules and the charging schedule of EVs, meanwhile subjecting to specific pickup and delivery time windows. In practice, the strict schedules put limitations on the operational potential of the fleet, hence, allowing for some time flexibility with respect to the desired schedule can significantly improve the overall operational performance of the entire system. To address this issue, the fleet operator can purchase such time flexibility from customers, and still obtain the larger overall revenues. From the perspective of the customers, if enough compensation is offered, it can be beneficial to provide time flexibility in exchange for a reduction of the delivery fee. In this context, to maximize the fleet operator's revenues, the fleet operator aims to seek a great balance between the monetary incentives offered to the customers and the benefits coming from the time flexibility provided by the customers. Note that with the adoption of an incentive design strategy into the EVs' routing problem (EVRP), there are two-fold benefits: (1) the customers will receive the delivery fee savings with the discount offered; (2) when the operator has an adjustable time window, more flexibility for the routing and charging schedule will be obtained resulting in the reduced operation cost. For this reason this paper investigates the EVRP with incentive-aware customers accounting for the decision of customers, and the monetary incentive mechanism design of the fleet operator.

\noindent\textit{Related work:} This paper organizes the relevant literature as the following two parts: (i) monetary incentive mechanism design, (ii) the EVRP. In order to leverage the power of monetary incentives, a few researchers focus on designing incentive schemes to optimize the EVs' charging problem and EVs' routing problem~\cite{narassimhan2018role,salazar2021urgency,lin2021electric,kucukoglu2021electric,basso2021electric,xiong2020integrated}. In order to find advanced transportation management solutions,  Xiong et al.~\cite{xiong2020integrated} developed an integrated and personalized traveler information scheme to stimulate an energy efficient mobility decision. Aiming to obtain a better routing and charging scheme, Diaz et al.~\cite{diaz2021integrated} proposed an incentive scheme, which enables the EV demand-aggregator have more room to devise a better operation scheme by determining the optimal charging time. However, these incentive schemes investigated in previous papers are solely used in either routing or charging applications, and there is no research publication on the joint optimization of the EV routing decisions, charging schedules, and incentive mechanism design problem.


Due to the local zero-carbon emissions of EVs, the EVRP has recently drawn great attention from the academic community~\cite{basso2021electric,yao2021joint,boewing2020vehicle,luke2021joint,lin2021electric,ferro2020optimal}.
In order to cope with the computational complexity stemming from the routing and charging problem of EVs, which is essentially a mixed-integer nonlinear programming problem, we proposed a computational-efficient method based on convex optimization~\cite{yao2021joint}. In the similar way, Boewing et al. devise a solution algorithm with polynomial time guarantees to obtain  the optimal vehicle routing and charging schedules of a fleet of electric autonomous vehicles offering mobility-on-demand service~\cite{boewing2020vehicle}. However, these papers do not incorporate the underlying interaction between the customers and the fleet operator, which can greatly enhance the performance of the entire system by reducing the operation cost for the fleet operator and helping the customers save more delivery fee. Given the time-varying characteristics of electricity pricing, there are some works focusing on the electric vehicle routing and charging problem with a time-of-use electricity price\cite{lin2021electric,ferro2020optimal}. To handle the EV routing problem with time windows considering time-variant electricity prices,  Lin et al.\cite{lin2021electric} proposed the Lagrangian relaxation approach and the tabu search method to obtain near optimal solutions. However, an impractical assumption that EV cannot start the charging process at a recharging station before the start of the following time period, is made in this paper. For this reason we aim to devise a mathematical model incorporating a spatio-temporal charging price without such an assumption.


In this paper, the mathematical formulation of an EVs' fleet operator is extended to incorporate the spatio-temporal electricity price at each charging station. Due to the hierarchical relation between fleet operator and customers, a bi-level optimization is employed which is equivalently reformulated as a single-level optimization problem. In practice, the customers (e.g., the e-commerce company, and the goods distributor) are prone to obtain the transportation service from the delivery companies (fleet operator) like FedEx and UPS, which can provide transportation as a service. Thus, considering the aforementioned business model and the privacy-preserving requirement of a fleet operator and customers, it is reasonable to decompose the resultant single-level optimization problem to subproblems which are individually solved by the freight transportation operator and the customers. To this end we proceed with the strengthened generalized Benders decomposition (SGBD) method. 


\noindent\textit{Contributions:} 
\begin{itemize}
    \item \textcolor{black}{ We accommodate the spatio-temporal electricity price and propose a detailed mathematical formulation for the fleet operator. This formulation minimizes the operation cost including charging cost, traveling time, charging time, EV usage cost, and delivery fee discount. }
    
    \item  We present a bi-level optimization model where the fleet operator aims to solve the EVRP with incentive-aware customers and determines the monetary incentives offered to the customers, which are also dependent on the time flexibility allowed by the customers.
    
    \item  To further cope with the nonlinear terms in the operator problem of the bi-level model, we propose a set of exact transformation techniques consisting of the KKT optimality condition, the Big-M linearization method, and the zero duality gap of the convex optimization problem to transform the proposed bi-level model to a single-level mixed-integer programming (MIP) problem.
    
    \item  To address the inherent complexities of the single-level MIP problem, and the privacy-preserving requirement of the customers and the fleet operator, we devise an SGBD which combines the merits of the generalized Benders decomposition (GBD) method and the Lagrangian dual decomposition method. We prove that the cuts from SGBD are much tighter than the cuts of GBD.


\end{itemize}



A preliminary version of this paper was accepted by the 2022 American Control Conference~\cite{yao2022incentive}. In this revised and extended version, we propose a detailed mathematical model including spatio-temporal electricity prices. Furthermore, leveraging the hierarchical structure of the customer problem and the fleet operator problem, we devise a decomposition-based algorithm combining the benefits of the generalized Benders decomposition and the Lagrangian dual decomposition method which iteratively solves the resulting single-level problem. Finally, additional numerical results regarding algorithmic performance are presented and discussed.

\noindent\textit{Organization:} The paper unfolds as follows. The system models and the mathematical model of the bi-level EVRP with incentive-aware customers and spatio-temporal electricity price are elaborated in Section~\ref{sec:mathematical}. In Section~\ref{sec:reformulation}, an exact reformulation technique is devised to reformulate the bi-level optimization problem as a single-level MIP. To achieve a decomposed implementation, we devise a decomposition algorithm with stronger valid cuts  in Section~\ref{sec:distributed}. 
Extensive numerical simulations are carried out in Section~\ref{sec:simus} to demonstrate the superiority of the proposed bi-level mathematical model. In Section~\ref{sec:conclusion}, we draw the conclusions of this paper. Note that the terms "requests" and "customers" are used interchangeably.


\section{Mathematical Formulation}\label{sec:mathematical}
To characterize a practical scenario of the EVRP with incentive-aware customers, we devise a bi-level optimization model to capture the interaction between the customers and the fleet operator. In this section we present the fleet operator model and the customers model as optimization problems minimizing the operational costs and the total inconvenience perceived by the customers, respectively. Finally, we combine these two optimization problems into a bi-level framework.

\subsection{Fleet Operator Model}

\subsubsection{The routing process constraint} We define the relation between electric vehicles and customers with a directed graph as in our previous work~\cite{yao2021joint}.
We denote the transportation network as a directed graph $G(\mathcal{V},\mathcal{E})$, where $\mathcal{V}=\{v_1,\mathcal{C},\mathcal{R},v_n\}$ comprises a start depot $v_1$, an end depot $v_n$, charging station nodes $\mathcal{C}$, and nodes representing the customers set $\mathcal{R}$.
The set $\mathcal{E}$ stands for the set of paths with $(i,j)\in\mathcal{E}$ representing a path from node $i$ to $j$. We define $d_{ij}$ as the travel distance between node $i\in \mathcal{V}$ and node $j\in \mathcal{V}$. Similarly, $T_{ij},e_{ij}$ represent travel time and energy consumption from node $i$ to $j$, respectively. The binary variable $x_{ij}^k$  shows whether vehicle $k$ is assigned to traverse path $(i,j)$. In order to capture the practical characteristics of the charging network and electric vehicles, the electric vehicle flow conservation constraints, battery energy constraints and visiting time constraints are clarified as follows: 

\par Here, all EVs are subject to the vehicle flow conservation constraint: An EV entering a customer node has to leave the same customer node, as well as  EVs starting at the start depot return back to the end depot after serving the requests. This is given by
  \begin{equation}\label{Cons_flow}
\begin{aligned}
\sum\limits_{j\in\mathcal{V}} & x^k_{i j}-\sum\limits_{j\in\mathcal{V}} x^k_{j i}=b_{i}, \quad \forall  i \in \mathcal{V}, k\in\mathcal{K},\quad\text{where}\quad b_{v_1}=1, b_{v_n}=-1, b_{i}=0.
\end{aligned}
\end{equation}
There is a constraint that each customer is served at most once by an EV as
\begin{equation}\label{Cons_visit}
    \sum\limits_{k\in\mathcal{K}}\sum\limits_{j\in\mathcal{V}} x^k_{ij} \leq 1,
    \quad\forall i\in\mathcal{R}.
\end{equation}
Because a few transportation customers may not be served by EVs if the serving cost is higher than the benefit occurring during the delivery service. Besides, visiting time constraints are characterized by
\begin{equation}\label{Cons_time}
	\begin{aligned}
		t_{j}   \geq &T_{ij}+g_i\cdot r^k_{i}+t_{i} -M(1-x^k_{ij}),\quad \forall i\in \mathcal{C}, j \in \mathcal{V}\setminus  v_1,k\in\mathcal{K},
		\end{aligned}
\end{equation}
\begin{equation}\label{Cons_time_nocharging}
	\begin{aligned}
		 t_{j}   \geq & T_{ij} + t_{i} -M(1-x^k_{ij}),\quad  \forall i\in \mathcal{R}, j \in \mathcal{V}\setminus  v_1,k\in\mathcal{K},	\end{aligned}
\end{equation}
in which \eqref{Cons_time} specifies that the arrival time at the following customer node $t_j$ should be longer than the total of the arrival time of previous customers $t_i$, charging time $g_i\cdot r_i^k$ in which $g_i$ and $r_i^k$ represent the time duration per unit kilowatt and the charging amount of EV $k$ at charging station $i$, respectively, and travel time $T_{ij}$. Constraint \eqref{Cons_time_nocharging} denotes that the arrival time of the following request node $t_j$ should be longer than the total of the arrival time of previous transportation requests $t_i$, and travel time $T_{ij}$. The constant $M$ of \eqref{Cons_time} and \eqref{Cons_time_nocharging} representing a constant with a large value (e.g., $M=10^6$), is introduced to avoid the occurrence of bilinear terms. The same techniques also apply to \eqref{Cons_soc} and \eqref{Cons_soc_nocharge}. In addition, each customer is constrained by their own time windows, which can be expressed as follows:
\begin{equation}\label{Cons_time-TW}
	\begin{aligned}
		t_{j}^{\mathrm{L}} \leq	t_{j} \leq  t_{j}^{\mathrm{L}} + \delta_{j},  \quad\forall  j \in \mathcal{R}.
	\end{aligned}
\end{equation}

\subsubsection{The charging process constraint}
In contrast to the frequently used spatially varying electricity price of the mathematical model in literature~\cite{chen2016optimal,yao2021joint}, we extend the mathematical model to incorporate the time-varying electricity price which is partly addressed by Lin et al.~\cite{lin2021electric} and Ferro at al.~\cite{ferro2020optimal}. However, they both have a similar but impractical assumption that the EVs need to recharge their batteries at the start of the time slot. Under the circumstance of the time-varying electricity price, since the start time of the charging process is dependent on the charging and routing schedules, it is hard to characterize the charging cost, which is the product of the charging energy and the time-varying electricity price. Extending the mathematical model of Ferro et al.~\cite{ferro2020optimal}, an additional continuous variable $r_i^k$ and two binary variables $B_{i\tau},B_{i\tau}^s$ are used to precisely characterize the relation between discrete time slots $\tau$ and the continuous arrival time of EVs $t_i$,  and to incorporate the spatio-temporal electricity cost.

\par The battery energy dynamics of the EVs at the charging station node is defined by
\begin{equation}\label{Cons_soc}
	\begin{aligned}
		-M (1 - &x^k_{ij})\leq -E^k_j+E^k_{i}+r_i^k-e_{i j} x^k_{i j} 
		 \leq M (1-x^k_{ij}),\quad \forall i\in\mathcal{C}, j\in\mathcal{V}\setminus  v_1,k\in\mathcal{K},
	\end{aligned}
\end{equation}
\noindent where $r_i^k$ is the amount charged at charging node $i$. Besides, the battery energy dynamics of the EVs at the customer node are characterized by
\begin{equation}\label{Cons_soc_nocharge}
	\begin{aligned}
		-M (1-x^k_{ij})&\leq -E^k_j+E^k_{i}-e_{i j} x^k_{i j} \leq M (1-x^k_{ij}), \quad\forall i\in\mathcal{R}, j\in\mathcal{V}\setminus  v_1,k\in\mathcal{K},
	\end{aligned}
\end{equation}
in which $e_{ij}=\phi d_{ij}$ denotes the amount of energy consumption from node $i$ to node $j$, and $\phi$ denotes the battery energy consumption per unit distance. As for constraint \eqref{Cons_soc} and \eqref{Cons_soc_nocharge}, the Big-M linearization approach is used to reformulate the original nonlinear constraints as the linear constraints. We put lower and upper bounds on the battery energy level of EV $k$, $E_i^k$   as
\begin{equation}\label{Cons_socbound}
    0\leq E^k_i\leq E^k_\mathrm{max}, \quad  i\in \mathcal{V},k\in\mathcal{K},
\end{equation}
where $E^k_\mathrm{max}$ is the battery capacity of each EV. The initial battery energy level of EVs is
\begin{equation}\label{Cons_initialsoc}
E^k_{v_1}=E^k_0 \quad\forall k\in\mathcal{K}.
\end{equation}
In addition, the time horizon is divided as a set of time slots $\Lambda = \{1,\cdots,\xi \}$ where $\xi$ represents the maximum time slots, with the constant length $\Delta\tau$. Two binary decision variables $B_{i\tau},B^\mathrm{s}_{i\tau}$ are introduced, defining whether charging station $i$ recharges in time slot $\tau$ or not, and whether charging station $i$ starts to recharge in time slot $\tau$ or not, respectively. Besides, the relation between $B_{i\tau}$ and $B^\mathrm{s}_{i\tau}$ is captured by
\begin{subequations}\label{Cons:upper}
    \begin{align}
  & B^\mathrm{s}_{i\tau} \geq B_{i\tau} -B_{i(\tau-1)},\quad\forall i\in\mathcal{C},\tau\in\Lambda, \label{binary_s}\\
  &     \sum\limits_{\tau\in\Lambda} B^\mathrm{s}_{i\tau}\leq 1,\quad\forall i\in\mathcal{C}, \label{binary_charging}
    \end{align}
  \end{subequations}
  which specifies that the maximum number of charging processes occurring during the entire time slot $\Lambda$ is one.  Note that in practice the charging station could provide the charging service several times over the entire horizon; in other words, during the whole horizon multiple charging processes could take place.
\begin{remark}
To accommodate the scenario in which multiple charging processes occur during the whole horizon, the charging nodes are duplicated with $N$ dummy nodes, which are connected by the corresponding charging nodes by setting the distance to $0$ where $N$ denotes the maximum number of charging processes.
\end{remark}

Considering that we can represent the charging time in two different forms, namely, \textit{(a)} as the product of charging amount and charging time per kilowatt hour $r_i^k g_i$ and \textit{(b)} the sum of all charging time slots over the entire time horizon $\sum\limits_{\tau\in\Lambda} B_{i\tau} \Delta\tau$,  the relation between $r_i^k g_i$ and $\sum\limits_{\tau\in\Lambda} B_{i\tau} \Delta\tau$ is given by

\begin{equation} \label{Cons:st:energycost}
	r^k_i g_i \leq \sum\limits_{\tau\in\Lambda} B_{i\tau} \Delta\tau, \quad\forall i \in \mathcal{C}.
\end{equation}
The relation between discrete time slot $\tau$ and continuous arrival time of EVs $t_i$ is characterized by
  \begin{subequations}\label{Cons:lower}
    \begin{align}
  &  \sum\limits_{\tau \in \Lambda} B^s_{i \tau} \tau \Delta \tau \leq t_{i}, \quad \forall i \in \mathcal{C} \\
  & t_{i} \leq \sum\limits_{\tau \in \Lambda} B^s_{i \tau}(\tau+1)\Delta \tau, \quad \forall i \in \mathcal{C}. 
  \end{align}
  \end{subequations}

\subsubsection{The objective function}
\par With the time flexibility $\delta_j$ allowed by the customers, the fleet operator optimizes the operational cost comprising of 
\begin{itemize}
    \item the charging cost $\sum\limits_{i\in\mathcal{R}}\sum\limits_{\tau\in\Lambda}  \frac{ p_{i\tau} B_{i\tau}\Delta\tau}{g_i}$, 
    
    \item  the EVs usage cost and delivery revenue
$\sum\limits_{k\in \mathcal{K}}\sum\limits_{i\in\mathcal{V}} \sum\limits_{j\in\mathcal{V}}  c_{i}    x^k_{i j}$, 

 \item  the travel time $\sum\limits_{k\in \mathcal{K}}\sum\limits_{i\in\mathcal{V}} \sum\limits_{j\in\mathcal{V}}    \omega_\mathrm{T} T_{ij}  x^k_{i j}$, 
 
 \item  the charging time $\sum\limits_{k\in \mathcal{K}}\sum\limits_{i\in\mathcal{V}} \sum\limits_{j\in\mathcal{V}}   \omega_\mathrm{T}  r^k_{i} g_{i}  x^k_{i j}$, 
 
 \item  the total discount $ \sum\limits_{j\in \mathcal{R}} q_j \delta^*_j \sum\limits_{i\in \mathcal{V}}\sum\limits_{k\in\mathcal{K}} x^k_{ij}$,
\end{itemize}
which stands for the cost paid for the customers that provide a flexible delivery time window. Here $\omega_\mathrm{T}$ denotes the value of time. Besides, $c_i$ represents a unified cost vector capturing both the vehicle usage fee $c_{v}$ and the negative delivery revenue $D_i$\footnote{Since the fleet operator would like to minimize the operational cost, we set the delivery revenue $D_i$ of the customers as a negative value.}:

\begin{equation}\label{Cons:para}
	c_i=\left\{
	\begin{array}{cc}
		{D_i,} & {\text { if } i\in\mathcal{R}} \\ 
		{c_{v},} & {\text { if } i=v_1}.
	\end{array}
	\right.
\end{equation}

To summarize, the fleet operator problem with incentive mechanism design is stated as follows:
\begin{problem}[Operator model]\label{Model_operator}
     \begin{equation}\notag\label{obj}
\begin{aligned}
\min\limits_{x^k_{i j},B_{i\tau},B^\mathrm{s}_{i\tau} \in\mathbb{B}, r_i^k, q_j , t_j \in\mathbb{R}}  & \sum\limits_{i\in\mathcal{R}}\sum\limits_{\tau\in\Lambda}  \frac{ p_{i\tau} B_{i\tau}\Delta\tau}{g_i} +\sum\limits_{k\in \mathcal{K}}\sum\limits_{i\in\mathcal{V}} \sum\limits_{j\in\mathcal{V}}   (c_{i}  + \omega_\mathrm{T} T_{ij}   +  \omega_\mathrm{T}  r^k_{i} g_{i}  ) x^k_{i j} +  \sum\limits_{j\in \mathcal{R}} q_j \delta^*_j \sum\limits_{i\in \mathcal{V}}\sum\limits_{k\in\mathcal{K}} x^k_{ij}\\ 
\centering &\text{s.t. }\eqref{Cons_flow}-\eqref{Cons:para}.
  \end{aligned}
\end{equation}
\end{problem}

\subsection{Customers Model}
We assume that customers aim to maximize the incentive received, and to minimize the inconvenience coming from their provided time flexibility $\delta_j$. Besides, we quantify the inconvenience of the customers with a monotonically increasing convex function $\mathcal{I}(\delta_j)$~\cite{yao2021vehicle}.
Specifically, given the discount price $q_j$ of customer $j$ determined by the fleet operator, we formalize the customer problem as follows:
\begin{problem}[Customer Model]\label{Model_general_customer}
\begin{equation}\label{P_customer} 
\begin{aligned}
\min_{\delta_j\in\mathbb{R}}\qquad   &  \mathcal{I}(\delta_j)-q_j \delta_j,    \\\text{s.t.}\quad &  0 \leq \delta_j \leq \Bar{\delta}_j
\end{aligned}
\end{equation}    
\end{problem}

Thereby, we define these non-negative dual variables $u_j,\sigma_j$ which are associated with the inequality constraints of Problem \ref{Model_general_customer}.

\subsection{Joint Scheduling Problem of the Fleet Operator and Customers}
With the fleet operator model and customers model, we formulate the routing and charging problem of EVs considering the monetary incentives as follows:
\begin{problem}[Joint routing and charging model of the fleet operator and customers]\label{Model_bilevel}
  \begin{equation}\notag
	\begin{aligned}
		\min\limits_{x^k_{i j},B_{i\tau},B^\mathrm{s}_{i\tau} \in\mathbb{B}, r_i^k, q_j , t_j \in\mathbb{R}} & \sum\limits_{i\in\mathcal{R}}\sum\limits_{\tau\in\Lambda}  \frac{ p_{i\tau} B_{i\tau}\Delta\tau}{g_i}
		+ \sum\limits_{k\in \mathcal{K}}\sum\limits_{i\in\mathcal{V}} \sum\limits_{j\in\mathcal{V}}   (c_{i}   + \omega_\mathrm{T} T_{ij}   +  \omega_\mathrm{T}  r^k_{i} g_{i}  ) x^k_{i j} +  \sum\limits_{j\in \mathcal{R}} q_j \delta^*_j \sum\limits_{i\in \mathcal{V}}\sum\limits_{k\in\mathcal{K}} x^k_{ij}\\
\text{s.t. }\delta^*_j \in \arg\min_{\delta_j}  \Big\{ &\mathcal{I}(\delta_j)-q_j \delta_j ,   0\leq\delta_j \leq \Bar{\delta}_j \Big\}, \forall j \in \mathcal{R}, 
\\&\quad\eqref{Cons_flow}-\eqref{Cons:para}.
\end{aligned}
\end{equation}
\end{problem}

Regarding the inherent difficulties raised by this bi-level optimization problem, there are a few studies attempting to solve the bi-level optimization problem\cite{lozano2017value,sinha2017review,allende2013solving}. As an alternative to solving the bi-level optimization problem\cite{cheng2019low} with a proposed bisection-based iterative approach, we propose to equivalently reformulate the bi-level optimization problem as a single-level optimization problem. 

\allowdisplaybreaks

\section{Equivalent Single-level Reformulation Approach}\label{sec:reformulation}
Due to the inherent NP-hardness brought by the bi-level optimization, thus, we cannot obtain the global optimal solution within the polynomial computation time~\cite{dempe2015bilevel}. For this reason directly solving Problem~\ref{Model_bilevel} is unpractical. To address this issue, we employ the KKT optimality conditions of the customers optimization problem to accurately reformulate Problem~\ref{Model_bilevel} as a single-level mixed-integer nonlinear programming (MINLP) problem. In addition, to simplify such a complicated MINLP, we propose an accurate linearization method, 
which uses the Big-M linearization method, and the strong duality of the convex optimization problem.

\subsection{Reformulation of Customer Problem}
In order to handle the NP-hardness of the proposed bi-level model, by leveraging the strong duality of the customer problem of Problem~\ref{Model_bilevel}, we represent the customer problem as a convex problem via its KKT optimality conditions:

\begin{subequations}\label{KKT}
\begin{align}
& \nabla\mathcal{I}(\delta_j) -q_j -\sigma_j +u_j =0,\forall j\in\mathcal{R} \label{KKT_stationary_delta}\\
& 0\leq u_j \perp (\Bar{\delta}_j-\delta_j) \geq 0, \forall j\in\mathcal{R}\label{KKT_complementary_u}\\
 &0\leq \sigma_j \perp \delta_j \geq 0,\forall j\in\mathcal{R}\label{KKT_complementary_sigma}
 \end{align}
\end{subequations}


\noindent where $ \perp $ is the complementarity operator. The stationarity conditions are specified in \eqref{KKT_stationary_delta}, whilst the primal feasibility constraint, the dual feasibility constraint, and complementary condition are given by~\eqref{KKT_complementary_u} and \eqref{KKT_complementary_sigma}.


Note that there is still a nonlinear term $q_j\delta_j$ in the objective function of Problem~\ref{Model_bilevel}, and in the complementarity constraints \eqref{KKT_complementary_u} and \eqref{KKT_complementary_sigma}. In order to cope with the computational burden brought by these nonlinear terms, we devise an exact linearization approach to represent such nonlinear terms in a computationally-efficient manner.

\subsection{Accurate Linearization Method of the Nonlinear Terms}
In order to further reduce the computational complexity of Problem~\ref{Model_bilevel}, an equivalent linearization method consisting of \text{(i)} introducing auxiliary binary variables to linearize the nonlinear complementary constraint, and \text{(ii)} exploiting the strong duality of the customer problem~\cite{wei2014energy} to accurately linearize the nonlinear term $q_j\delta_j$, is devised.

\subsubsection{Linearized Complementary Constraints} To cope with the nonlinear and nonconvex complementary constraints \eqref{KKT_complementary_u}, and \eqref{KKT_complementary_sigma},  rendering the resulting problem difficult to solve, we devise an exact linearization approach which linearizes these hard constraints by introducing auxiliary binary variables $\psi^1_j,\psi^2_j$ and a constant $M$ with the sufficiently large value, resulting in the subsequent disjunctive constraints~\cite{fortuny1981representation}:

\begin{subequations}\label{Optimality Con linearized}
\begin{align}
&
\left.
\begin{aligned}
&0\leq\Bar{\delta}_j-\delta_j \leq M \psi^{1}_j\\
&0\leq u_j \leq M (1- \psi^{1}_j)
\end{aligned}
\right\}  \forall j\in\mathcal{R}
\\&
\left.
\begin{aligned}
	&0\leq \delta_j \leq M \psi^{2}_j\\
	&0\leq \sigma_j \leq M (1- \psi^{2}_j)
\end{aligned}
\right\}  \forall j\in\mathcal{R}
\end{align}
\end{subequations}

\subsubsection{Linearized Objective Function}

Since there is still a nonlinear term $q_j\delta_j$, which is the product of two continuous variables, in the objective function of  Problem \ref{Model_bilevel}.  With the customer model in Problem \ref{Model_general_customer}, in which $\mathcal{I}(\delta_j)$ is a convex function, the Lagrangian function can be obtained as shown below:
\begin{equation}\notag
\begin{aligned}
L(\delta_j,u_j,\sigma_j) = &\mathcal{I}(\delta_j)-q_j \delta_j - \sigma_j \delta_j  +u_j (\delta_j-\Bar{\delta}_j)\\& = \mathcal{I}(\delta_j)  +(u_j -q_j - \sigma_j  )\delta_j- u_j \Bar{\delta}_j
\end{aligned}
\end{equation}

The dual function can be obtained as follows:

\begin{equation}\notag
\begin{aligned}
g(u_j,\sigma_j) &= \inf_{\delta_j}  \mathcal{I}(\delta_j)  +(u_j -q_j - \sigma_j  )\delta_j- u_j \Bar{\delta}_j \\
& = - u_j \Bar{\delta}_j+ \inf_{\delta_j}  \phi(\delta_j) = - u_j \Bar{\delta}_j+ \phi^*(\delta^*_j)
\end{aligned}
\end{equation}

\noindent For convenience, we denote $\phi(\delta_j) =   \mathcal{I}(\delta_j)  +(u_j -q_j - \sigma_j  )\delta_j $.

According to the strong duality of convex optimization problems, there is no duality gap between the objective value of lower-level Problem~\ref{Model_general_customer} and the objective value of its dual problem. Then the complicating nonlinear term $q_j\delta_j$ can be exactly linearized as follows:

\begin{equation}\notag
\begin{aligned}
 \mathcal{I}(\delta_j)-q_j \delta_j  &= - u_j \Bar{\delta}_j+ \phi^*(\delta^*_j)\\
     q_j \delta_j&= \mathcal{I}(\delta_j)+ u_j \Bar{\delta}_j- \phi^*(\delta^*_j)
\end{aligned}
\end{equation}

In addition, to handle the bilinear terms $\omega_\mathrm{T} r_i^k g_i x^k_{ij}$  and $(\mathcal{I}(\delta_j)+ u_j \Bar{\delta}_j- \phi^*(\delta^*_j))x^k_{ij}$ in the objective function of Problem~\ref{Model_bilevel}, the continuous auxiliary variables $\eta^1_{ijk},\eta^2_{j}$ are introduced as
\begin{subequations}\label{obj_linearized}
	\begin{align}
		&\eta^1_{ijk}\geq\omega_\mathrm{T}r^k_ig_i-M(1-x^k_{ij}),\forall i\in\mathcal{V},j\in\mathcal{V},k\in\mathcal{K}\\
		& \begin{aligned}
     	\eta^2_{j}\geq  \mathcal{I}(\delta_j)+ u_j \Bar{\delta}_j- \phi^*(\delta^*_j)&-M(1-\sum\limits_{k\in\mathcal{K}}\sum\limits_{i\in\mathcal{J}} x^k_{ij}),\forall j\in\mathcal{R}.
		\end{aligned}
	\end{align}
\end{subequations}

\subsubsection{Equivalent Reformulation of the bi-level optimization model of operator and customers}
We denote all variables by 
$ \mathcal{X}=\big\{\{x^k_{i j},\eta^1_{ijk}\}_{i,j\in\mathcal{V}}^{k\in\mathcal{K}}, \{q_j,\delta_j, u_j,\sigma_j, \psi^1_j,\psi^{2}_j,\eta^2_{j} \}_{j\in\mathcal{V}}, \{r_i^k,E_i^k\}_{i\in\mathcal{V}}^{k\in\mathcal{K}},$ $ \{B_{i\tau},B^\mathrm{s}_{i\tau} \}_{i\in\mathcal{C}}^{\tau\in\Lambda }  \big\}$ for the sake of convenience. With the proposed equivalent linearization methods, we can accurately reformulate Problem~\ref{Model_bilevel} as a Problem~\ref{Model_single_level_MILP}:
  \begin{problem}[Linearized single-level optimization problem]\label{Model_single_level_MILP}
      \begin{equation}\notag
\begin{aligned}
\min\limits_{\mathcal{X}} \sum\limits_{k\in \mathcal{K}}\sum\limits_{i\in\mathcal{V}} \sum\limits_{j\in\mathcal{V}} & \Big( \eta^1_{ijk}+ \eta^2_{j}  + ( \omega_\mathrm{T} T_{ij} + c_{i}) x^k_{i j} + \sum\limits_{i\in\mathcal{R}}\sum\limits_{\tau\in\Lambda}  \frac{ p_{i\tau} B_{i\tau}\Delta\tau}{g_i}  \Big)
 \end{aligned}
\end{equation}
$$\textit{s.t.}\quad \eqref{Cons_flow}-\eqref{Cons:para},\eqref{KKT_stationary_delta}, \eqref{Optimality Con linearized},\eqref{obj_linearized}. $$
  \end{problem}

Commercial solvers (e.g., Gurobi, and CPLEX) can be used to solve Problem~\ref{Model_single_level_MILP} directly. However, considering that the privacy-preserving requirement of both the fleet operator and the customers, we propose a decomposition-based algorithm to decompose the single-level optimization problem into a master problem and subproblem, which are solved in an iterative fashion. The master problem is solved by the fleet operator, and the subproblem is solved by a cloud operator of all customers. We assume that a cloud operator can collect all information regarding the solution of the subproblem. To book a ride-hailing service online, in reality, the customers need to upload their travel information to the online ride-hailing platform, like Uber. In return, the online ride-hailing service providers offer customers multiple travel options. Then the assumption aforementioned is reasonable.


\section{Optimization Framework with Decomposition}\label{sec:distributed}
Considering the independence of the EVs fleet operator and the customers, it is impractical to directly solve the resultant single-level Problem \ref{Model_single_level_MILP} in a centralized fashion. In the community of optimization, enormous efforts are put into designing efficient algorithms to solve mixed-integer programming problems~\cite{lee2011mixed,dakin1965tree,bodur2017strengthened,rahmaniani2020benders,bodur2017mixed,baena2020stabilized}. Due to the hierarchical relation between the fleet operator and the customers, we propose an optimization algorithm based on the SGBD, which combines the complementary merits of the Lagrangian dual decomposition and the GBD. More specifically, by exploiting the merits of the Lagrangian dual decomposition method, we devise an  MIP subproblem. With the newly formulated subproblem, valid Benders cuts are generated that are tighter than the classic Benders cuts derived from the GBD. Besides, the tightness is also rigorously proved in the subsection below.

\subsection{Generalized Benders Decomposition Method}

For the sake of simplicity, the linearized single-level Problem \ref{Model_single_level_MILP} is presented in compact form as follows: 
\begin{subequations} 
\begin{align}
    \min\limits_{X_\mathrm{d},X_\mathrm{c}}\quad & c_\mathrm{d}^T X_\mathrm{d} + c_\mathrm{c}^T X_\mathrm{c} \notag\\
    s.t.\quad & A_\mathrm{d} X_\mathrm{d} + A_\mathrm{c} X_\mathrm{c} \leq b_\mathrm{a}  \label{cons:MILP}\\
              & B_\mathrm{d} X_\mathrm{d} + G(X_\mathrm{c}) \leq b_\mathrm{b}  \label{cons:MIP}\\
              & D_\mathrm{d} X_\mathrm{d}  \leq b_\mathrm{c}  \label{cons:IP} \\
              & X_\mathrm{d} \in \mathbb{B}^m,  X_\mathrm{c} \in\mathbb{R}^n, 
    \end{align}
\end{subequations}
where $X_\mathrm{c}=\big\{\{\eta^1_{ijk}\}_{i,j\in\mathcal{V}}^{k\in\mathcal{K}}, \{q_j,\delta_j, u_j,\sigma_j,\eta^2_{j} \}_{j\in\mathcal{V}}$,$ \{r_i^k,E_i^k\}_{i\in\mathcal{V}}^{k\in\mathcal{K}}\big\}$ and $X_\mathrm{d}=\big\{\{x^k_{i j}\}_{i,j\in\mathcal{V}}^{k\in\mathcal{K}}, \{ \psi^1_j,\psi^{2}_j\}_{j\in\mathcal{V}},\{B_{i\tau},B^\mathrm{s}_{i\tau} \}_{i\in\mathcal{C}}^{\tau\in\Lambda }  \big\}$
 represent the set of continuous variables and the set of discrete variables, respectively. Besides, $A_\mathrm{d},A_\mathrm{c},B_\mathrm{d},D_\mathrm{d}$ are the coefficient matrices, $b_\mathrm{a},b_\mathrm{b},b_\mathrm{c},c_\mathrm{d},c_\mathrm{c}$ are the coefficient vectors, and $n,m$ represent the dimension of continuous variables and binary variables, respectively. In addition, \eqref{cons:MILP}  denotes the linear constraint, \eqref{cons:MIP} denotes the nonlinear constraint, where $G(X_\mathrm{c})$ represents a nonlinear function, and \eqref{cons:IP} represents the integer related constraints. 


\subsubsection{The mathematical formulation of the Benders master problem}
By adopting an auxiliary variable $\Theta$ serving as a lower bound of the objective value of the subproblem, the Benders master problem (BMP) containing integer variables~$X_\mathrm{d}$  and the related constraints~\eqref{cons:IP} can be formulated as follows:

\begin{problem}[BMP]\label{BMP}
\begin{equation}\notag
\begin{aligned}
    \min\limits_{X_\mathrm{d},\Theta}\quad & c_\mathrm{d}^T X_\mathrm{d} + \Theta \\
    s.t.\quad & D_\mathrm{d} X_\mathrm{d}  \leq b_\mathrm{c} \\
              & X_\mathrm{d} \in \mathbb{B}^m. \\
    \end{aligned}
\end{equation}
\end{problem}

\subsubsection{The mathematical formulation of the Benders subproblem}
With the solution obtained from solving BMP $X_\mathrm{d}^*$ and the local copies of BMP variables $Z_\mathrm{d}$, the Benders subproblem (BSP) can be formalized thus:

\begin{problem}[BSP]\label{BSP}
\begin{subequations}
\begin{align}
    \min\limits_{X_\mathrm{c},Z_\mathrm{d}}\quad &  c_\mathrm{c}^T X_\mathrm{c} \notag\\
    s.t.\quad & A_\mathrm{d} Z_\mathrm{d} + A_\mathrm{c} X_\mathrm{c} \leq b_\mathrm{a} \\
              & B_\mathrm{d} Z_\mathrm{d} + G(X_\mathrm{c}) \leq b_\mathrm{b} \\
              & D_\mathrm{d} Z_\mathrm{d}  \leq b_\mathrm{c} \\
              &  Z_\mathrm{d} = X_\mathrm{d}^* \label{cons:integrality}\\
              & Z_\mathrm{d} \in \mathbb{R}^m,  X_\mathrm{c} \in\mathbb{R}^n.
    \end{align}
\end{subequations}
\end{problem}

If the solution of Problem~\ref{BSP} $( \bar{X}_\mathrm{c} , \bar{Z}_\mathrm{d})$  is feasible, the following Benders optimality cut can be generated:

\begin{equation}\label{cut:bsp}
    \Theta \geq c_\mathrm{c}^T \bar{X}_\mathrm{c} + \zeta^T(X_\mathrm{d}- \bar{Z}_\mathrm{d}),
\end{equation}
where $\zeta$ is the dual variable of \eqref{cons:integrality}.

As for the infeasible Problem \ref{BSP}, with two slack variables $S_\mathrm{a}, S_\mathrm{b}$, the feasibility subproblem (FSP) is formulated as follows: 
\begin{problem}[FSP]\label{FSP}
\begin{subequations}
\begin{align}
    \min\limits_{X_\mathrm{c},Z_\mathrm{d},S_\mathrm{a},S_\mathrm{b}}\quad & \mathds{1}^T S_\mathrm{a} + \mathds{1}^T S_\mathrm{b}  \notag\\
    s.t.\qquad & A_\mathrm{d} Z_\mathrm{d} + A_\mathrm{c} X_\mathrm{c} \leq b_\mathrm{a} +S_\mathrm{a} \\
              & B_\mathrm{d} Z_\mathrm{d} + G(X_\mathrm{c}) \leq b_\mathrm{b} + S_\mathrm{b} \\
              & D_\mathrm{d} Z_\mathrm{d}  \leq b_\mathrm{c}  \\
              & Z_\mathrm{d} = X_\mathrm{d}^* \label{cons:integrality_f}\\
              & Z_\mathrm{d} \in \mathbb{R}^m,  X_\mathrm{c} \in\mathbb{R}^n.
    \end{align}
\end{subequations}
\end{problem}

With the solution of Problem~\ref{FSP} $(\bar{S}_\mathrm{a}, \bar{S}_\mathrm{b},  \bar{Z}_\mathrm{d})$, the Benders feasibility cut can be generated as follows:
\begin{equation}\label{cut:fsp}
	0 \geq \mathds{1}^T \bar{S}_\mathrm{a} + \mathds{1}^T \bar{S}_\mathrm{b}  + \Upsilon^T(X_\mathrm{d}- \bar{Z}_\mathrm{d}),
\end{equation}
where $\Upsilon$ is the dual variable of \eqref{cons:integrality_f}.

\subsection{Strengthened Generalized Benders Decomposition Algorithm}
Considering that directly applying the GBD is a time-consuming task\cite{lee2020accelerating}, we devise a novel decomposition method, namely, the SGBD, which iteratively strengthens the master problem with tight optimality and feasibility cut by leveraging the power of the Lagrangian dual decomposition method.

\subsubsection{Strengthened optimality cut}  There are a few papers developing the Benders optimality and feasibility cut generated at fractional nodes of the search tree during the early stage of the search process, which helps to obtain a better lower bound\cite{adulyasak2015benders,bodur2017strengthened,gendron2016branch}. In order to derive the strengthened optimality cut, by introducing the dual variable $\zeta$ for \eqref{cons:integrality}, the Lagrangian dual problem of \eqref{BSP} can be formulated as follows:
\begin{subequations}\notag
\begin{align}
\max\limits_{\zeta\in\mathbb{R}^m}    \min\limits_{X_\mathrm{c},Z_\mathrm{d}}\quad &  c_\mathrm{c}^T X_\mathrm{c}- \zeta^T(Z_\mathrm{d}- X_\mathrm{d}^*)  \notag\\
    s.t.\quad & A_\mathrm{d} Z_\mathrm{d} + A_\mathrm{c} X_\mathrm{c} \leq b_\mathrm{a} \\
              & B_\mathrm{d} Z_\mathrm{d} + G(X_\mathrm{c}) \leq b_\mathrm{b} \\
              & D_\mathrm{d} Z_\mathrm{d}  \leq b_\mathrm{c} \\
              & Z_\mathrm{d} \in \mathbb{R}^m,  X_\mathrm{c} \in\mathbb{R}^n.
    \end{align}
\end{subequations}

The strengthened Benders optimality cut can be generated by  enforcing $Z_\mathrm{d}\in\mathbb{B}^m$ within the Lagrangian dual problem of \eqref{BSP}, as illustrated in Lemma~\ref{cut:sopt}.

\begin{lemma}\label{cut:sopt}
Given the linear relaxation solution of BMP $X^*_\mathrm{d}\in\mathbb{R}^m$, and the dual variable $\zeta\in\mathbb{R}^m$ from solving BSP, we can formulate the following mixed-integer programming subproblem
\begin{subequations}\label{P:ssp}
\begin{align}
\min\limits_{X_\mathrm{c},Z_\mathrm{d}}\quad &  c_\mathrm{c}^T X_\mathrm{c} - \zeta^T(Z_\mathrm{d}- X_\mathrm{d}^*) \notag\\
    s.t.\quad & A_\mathrm{d} Z_\mathrm{d} + A_\mathrm{c} X_\mathrm{c} \leq b_\mathrm{a} \label{ssp:a}\\
              & B_\mathrm{d} Z_\mathrm{d} + G(X_\mathrm{c}) \leq b_\mathrm{b} \label{ssp:b}\\
              & D_\mathrm{d} Z_\mathrm{d}  \leq b_\mathrm{c} \label{ssp:c}\\
              &  Z_\mathrm{d} \in \mathbb{B}^m,  X_\mathrm{c} \in\mathbb{R}^n.
    \end{align}
\end{subequations}
to derive the strengthened Benders optimality cut \eqref{cuts:ssp} with the solution $\bar{Z}_\mathrm{d},\bar{X}_\mathrm{c}$ from solving \eqref{P:ssp}.
\begin{equation}\label{cuts:ssp}
    \Theta \geq c_\mathrm{c}^T \bar{X}_\mathrm{c} + \zeta^T(X_\mathrm{d}- \bar{Z}_\mathrm{d}).
\end{equation}
\end{lemma}
\begin{proof}
The proof can be found in Appendix~\ref{proof:lemma}.
\end{proof}

Compared with the optimality cut generated by the GBD, the tightness of the strengthened Benders optimality cut is quantified and rigorously proved in Theorem \ref{cut:gap}. 

\begin{theorem}\label{cut:gap}
Given the dual multipliers $\bar{\zeta}$ of \eqref{cons:integrality}  from solving Problem \eqref{BSP},  the strengthened Benders optimality cut \eqref{cuts:ssp} is parallel to the generalized Benders optimality cut \eqref{cut:bsp} and at least $\Xi \geq 0$ units tighter, where
 \begin{equation}\notag
     \Xi =  \min\limits_{\mathcal{X}\cup Z_\mathrm{d}\in\mathbb{B}^m  }(c_\mathrm{c}^T X_\mathrm{c} - \bar{\zeta}^T Z_\mathrm{d}  )   - \min\limits_{\mathcal{X}\cup Z_\mathrm{d}\in\mathbb{R}^m  }(c_\mathrm{c}^T X_\mathrm{c} - \bar{\zeta}^T Z_\mathrm{d}  ).    
 \end{equation} 
\end{theorem}
\begin{proof}
The proof is delivered in Appendix~\ref{proof:theorem}.
\end{proof}

With the tightly strengthened Benders cut of the SGBD, a larger area will be cut compared to the cut of the GBD. Thus, fewer iterations are achieved by the SGBD.

\subsubsection{Strengthened feasibility cut} As for the infeasible strengthened Benders subproblem, the strengthened feasibility subproblem (SFSP) of SGBD is formulated with three additional non-negative slack variables $(S_\mathrm{a},S_\mathrm{b},S_\mathrm{c})$ as follows:
\begin{problem}[SFSP]\label{P:SFSP}
\begin{subequations}
\begin{align}
    \min\limits_{X_\mathrm{c},Z_\mathrm{d},S_\mathrm{a},S_\mathrm{b},S_\mathrm{c}}\quad & \mathds{1}^T S_\mathrm{a} + \mathds{1}^T S_\mathrm{b}+ \mathds{1}^T S_\mathrm{c} - \Upsilon^T (Z_\mathrm{d}- X_\mathrm{d}^*)  \notag\\
    s.t.\quad & A_\mathrm{d} Z_\mathrm{d} + A_\mathrm{c} X_\mathrm{c} \leq b_\mathrm{a} +S_\mathrm{a} \\
              & B_\mathrm{d} Z_\mathrm{d} + G(X_\mathrm{c}) \leq b_\mathrm{b} + S_\mathrm{b} \\
              & D_\mathrm{d} Z_\mathrm{d}  \leq b_\mathrm{c}  +S_\mathrm{c} \\
              & Z_\mathrm{d} \in \mathbb{B}^m,  X_\mathrm{c} \in\mathbb{R}^n.
    \end{align}
\end{subequations}
\end{problem}

 With the solution of Problem~\ref{P:SFSP} $(\bar{S}_\mathrm{a}, \bar{S}_\mathrm{b}, \bar{S}_\mathrm{c}, \bar{Z}_\mathrm{d})$, the strengthened Benders feasibility cut can be generated as follows:
\begin{equation}\label{cuts:fsp}
	0 \geq \mathds{1}^T \bar{S}_\mathrm{a} + \mathds{1}^T \bar{S}_\mathrm{b}  + \mathds{1}^T \bar{S}_\mathrm{c} + \Upsilon^T(X_\mathrm{d}- \bar{Z}_\mathrm{d}).
\end{equation}

With these strengthened Benders cuts obtained~\eqref{cuts:ssp}, \eqref{cuts:fsp}, compared with the GBD, the performance of the proposed SGBD is theoretically improved, which will be evaluated below, in the numerical simulation results. 

\subsection{Implementation Details}

In this subsection, the implementation details are clarified as follows:

\begin{itemize}
    \item  Fractional solution: To quickly derive valid cuts, we first solve the LP relaxation of the BMP (Problem \ref{BMP}) with classical cuts;
    \item Strengthened cut generation: The strengthened Benders cuts \eqref{cuts:ssp}, \eqref{cuts:fsp} are generated in the following steps: (a) the value of dual multipliers  can be obtained by solving the subproblem \eqref{P:ssp} or the feasibility subproblem \eqref{P:SFSP}; (b) with the dual multipliers, strengthened Benders cuts \eqref{cuts:ssp}, \eqref{cuts:fsp}  can be generated.
\end{itemize}

 Besides, the master problem is solved by the fleet operator, and the strengthened subproblem is solved by a cloud operator of the customers. Thanks to the power of the strengthened optimality and feasibility cut, compared with the GBD, the proposed SGBD holds a faster convergence rate over the whole solution process.


\section{Simulation Experiments}\label{sec:simus}
In this section we demonstrate the effectiveness of the proposed mathematical model with an example of package delivery. Besides, we provide the details of the numerical data used.
\subsection{Illustrative Parameter Settings}\label{sec:pars}
The VRP-REP data for Belgium\footnote{\url{http://www.vrp-rep.org/datasets/item/2017-0001.html}}, which consists of a 1000-node map with their geographic coordinates, out of which we choose a finite number of customers in a random fashion, is used in the following numerical experiments. As for the specifications of EVs, the battery capacity and recharging power are set as \SI{90}{\kWh}, \SI{150}{\kW} (fast charging), and \SI{22}{\kW} (slow charging)\footnote{\url{https://www.tesla.com/}}, respectively.  In addition, the energy consumption per kilometer and the average speed of EVs are set as \SI[per-mode=symbol]{0.24}{\kWh\per\km} and \SI[per-mode=symbol]{60}{\km\per\hour}, respectively. Considering the 24-hour planning horizon, there are 288 time slots with 5-minute time intervals. According to the real time electricity price of PJM wholesale market\footnote{\url{https://hourlypricing.comed.com/live-prices/five-minute-prices/}}, we obtain the electricity price of 288 time slots during 5 minutes for the numerical simulation. We set the usage cost of an EV $c_{v}$ as \SI[per-mode=symbol]{99}{\$}. In order to capture a realistic scenario, the positive revenue gained from serving a request, is sampled from a truncated normal distribution with the  mean and the standard deviation set as  \SI[per-mode=symbol]{9.05}{\$} and \SI{5}{\$}\footnote{\url{https://www.fedex.com/en-us/shipping/one-rate.html}}, respectively. For the sake of convenience, we model the inconvenience function as a convex two-segment piecewise affine function with $\gamma_{1,2} = \{ 0,1.5\}\si{\$\per\hour} $ and  $\chi_{1,2} = \{0.01,-0.01\}\si{\$} $ for every customer. Furthermore, we can use the different inconvenience function for each customer to incorporate the heterogeneous case (i.e., every customer has the different time sensitivity $\gamma_i$). The simulation experiments are run with 50 times, and we showcase the results with average values in the subsequent figures. We solve the optimization problems on a workstation with an Intel Core i9-10980XE processor consisting of 36 CPUs of \SI{3.00}{\giga\hertz} and \SI{64}{\giga\byte} of memory. 






  


\subsection{Performance Comparison}
To evaluate the benefit of introducing monetary incentives both for the fleet operator and the customers, with different values of the maximum time flexibility $\bar{\delta}_j$ (which we set to be the same value for all customers), and with different values of the time sensitivity $\gamma_2$, the optimal solution of Problem~\ref{Model_single_level_MILP}, and the solution of the incentive-free problem detailed in Appendix~\ref{app:noincentives}, are compared and analyzed.


\subsubsection{Benefits brought by the time-flexibility $\bar{\delta}$ }\label{Si:time Flexibility on operation cost} 
\begin{figure}
     \centering
     \begin{subfigure}[b]{0.49\textwidth}
         \centering
 \includegraphics[width=\linewidth]{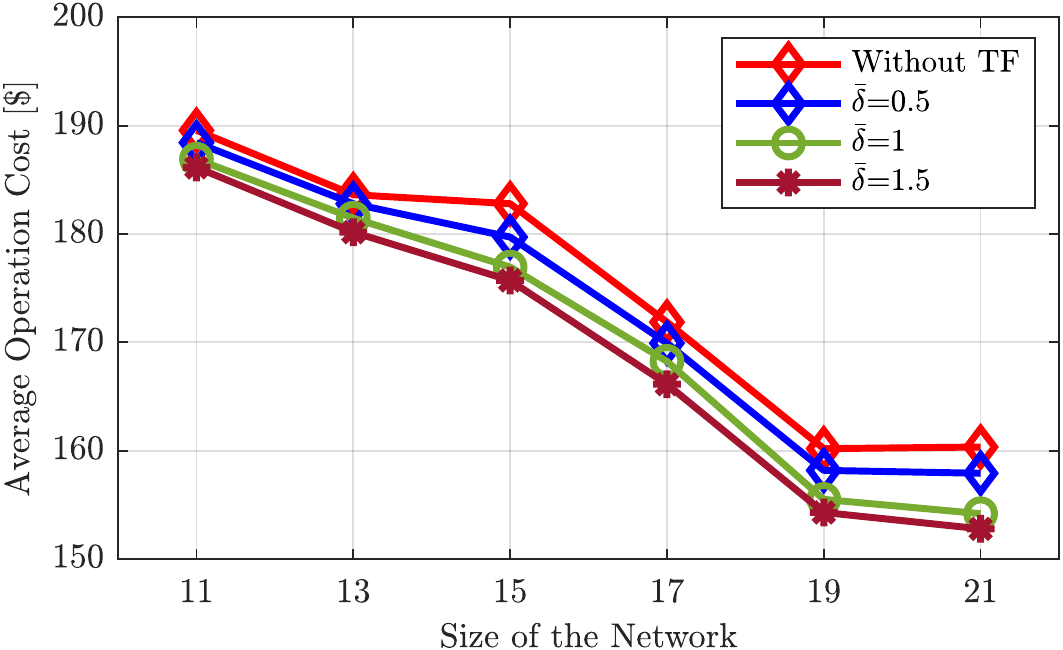}
     \caption{The average operation cost. }
     \label{OperationCostReduction_delta}
     \end{subfigure}
     \hfill
     \begin{subfigure}[b]{0.49\textwidth}
         \centering
\includegraphics[width=\linewidth]{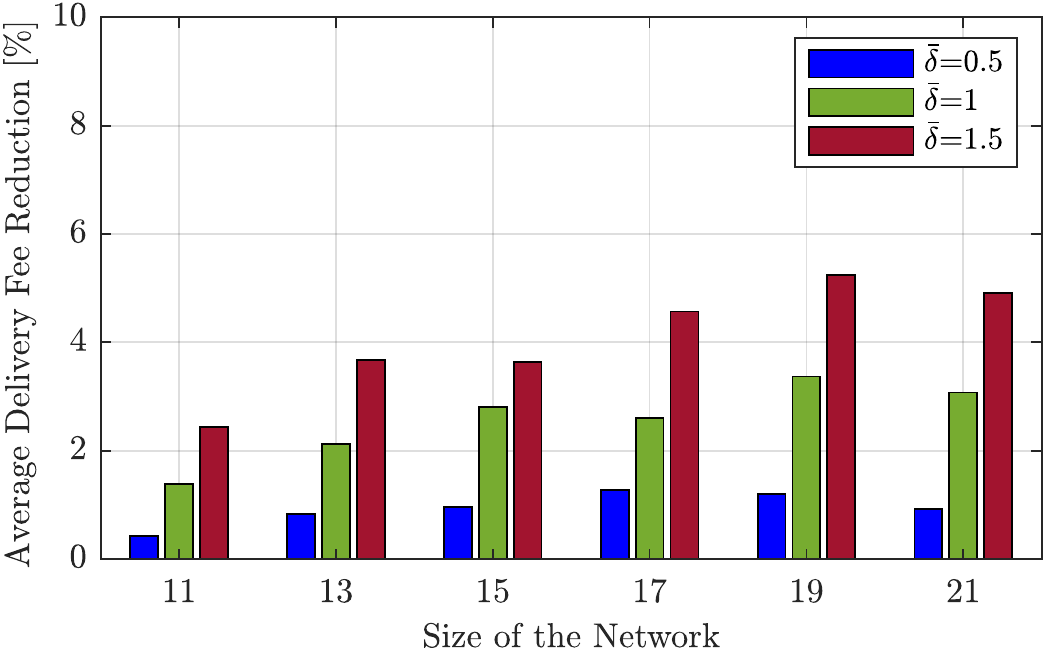}
 \caption{The average delivery fee savings. }
\label{DeliveryFeeReduction_delta}
     \end{subfigure}
     \caption{The impact of the different  $\bar{\delta}$ on the reduction of operational cost and the saving of delivery fee. }
\label{Time_flexibility}
\end{figure}

In order to assess the impact of different lengths of time flexibility on the reduction of the operation cost, we conduct 3 sets of experiments with the different value of $\bar{\delta}$,   $\bar{\delta}\in\{0.5, 1, 1.5\}\si{\hour}$, and $\gamma_2=$ \SI{1.5}{\$\per\hour}. As illustrated in Fig.~\ref{OperationCostReduction_delta}, the operation cost is decreased while the value of $\bar{\delta}$ increases from \SI{0.5}{\hour} to \SI{1.5}{\hour}, demonstrating that the operation cost can be further reduced as the customers provide a larger time flexibility.
In turn, a higher degree of time flexibility $\bar{\delta}$ will enable customers to receive a larger monetary incentive, ultimately leading to a larger delivery fee saving, as showcased in Fig.~\ref{DeliveryFeeReduction_delta}. Thereby, the maximum average operation cost reduction is over $5\%$, with the largest average delivery fee saving exceeding $5\%$.

\subsubsection{Benefits brought by the different time-sensitivity $\gamma_2$}\label{Si:time sensitivity on operation cost}  
\begin{figure}
     \centering
     \begin{subfigure}[b]{0.49\textwidth}
         \centering
 \includegraphics[width=\linewidth]{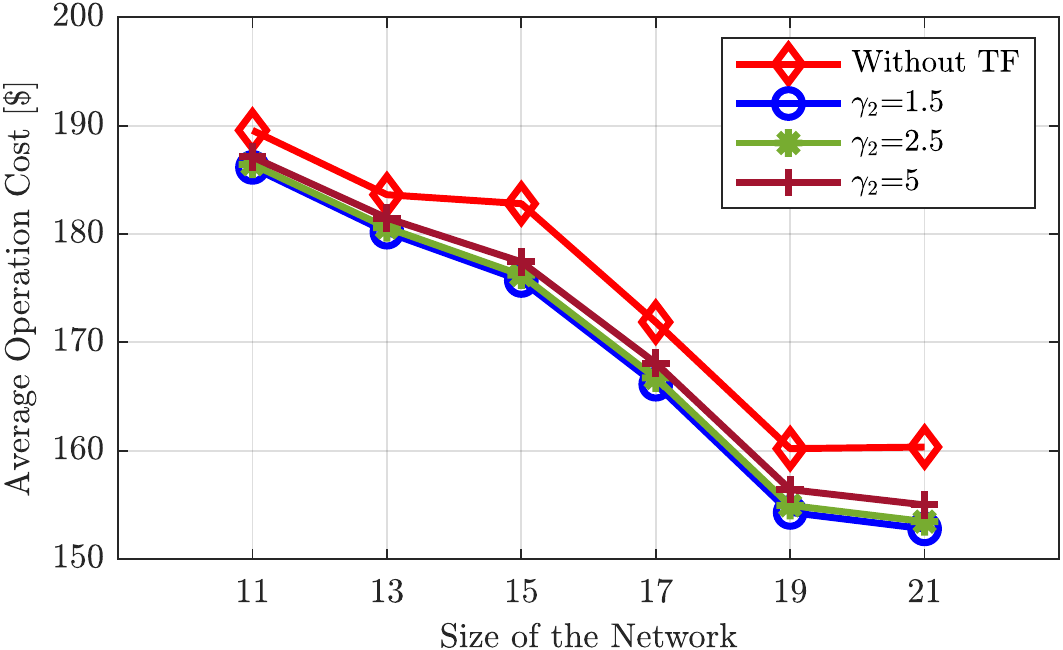}
     \caption{ The average operation cost.  }
       \label{OperationCostReduction_gamma}

     \end{subfigure}
     \hfill
     \begin{subfigure}[b]{0.49\textwidth}
         \centering
\includegraphics[width=\linewidth]{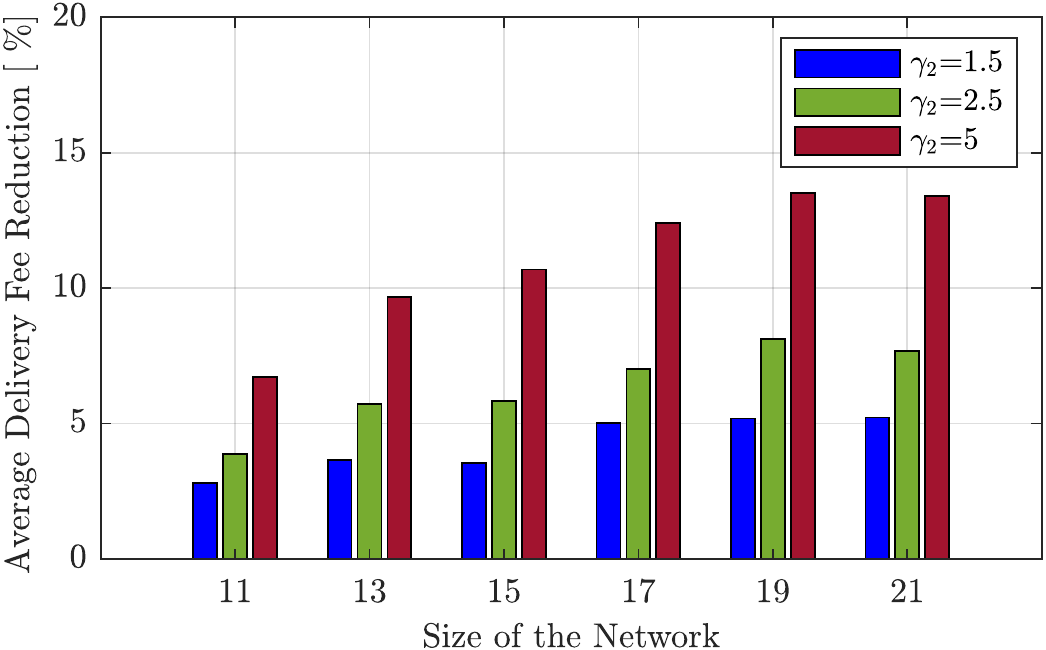}
     \caption{The average delivery fee savings.  }
     \label{DeliveryFeeReduction_gamma}
     \end{subfigure}
     \caption{The impact of different $\gamma_2$ on the reduction of operational cost and the saving of delivery fee. }
\end{figure}

We conduct three sets of simulation experiments with $\gamma_2 \in \{1.5,2.5,5\}\si{\$\per\hour}$ and $\bar{\delta}=$ $\SI{1.5}{\hour}$ to showcase the impact of the different degrees of time sensitivity of customers on the reduction of the operation cost.
As shown in Fig.~\ref{OperationCostReduction_gamma}, the value of the average operation cost becomes larger with the increasing value of $\gamma_2$.
Thus, as the customers pay more attention on the time flexibility, the fleet operator has to provide a higher monetary incentive in exchange for the time flexibility. Besides, we also conclude that even in the extreme case, it is always beneficial for the fleet operator to provide such incentives to acquire the flexibility of variable time windows. As illustrated in Fig.~\ref{DeliveryFeeReduction_gamma}, a larger value of $\gamma_2$ results in a larger average fee saving obtained, exceeding $13\%$. 



\subsubsection{Comparison between GBD and SGBD}
\begin{figure}[t]
\centering
\includegraphics[width=0.6\linewidth]{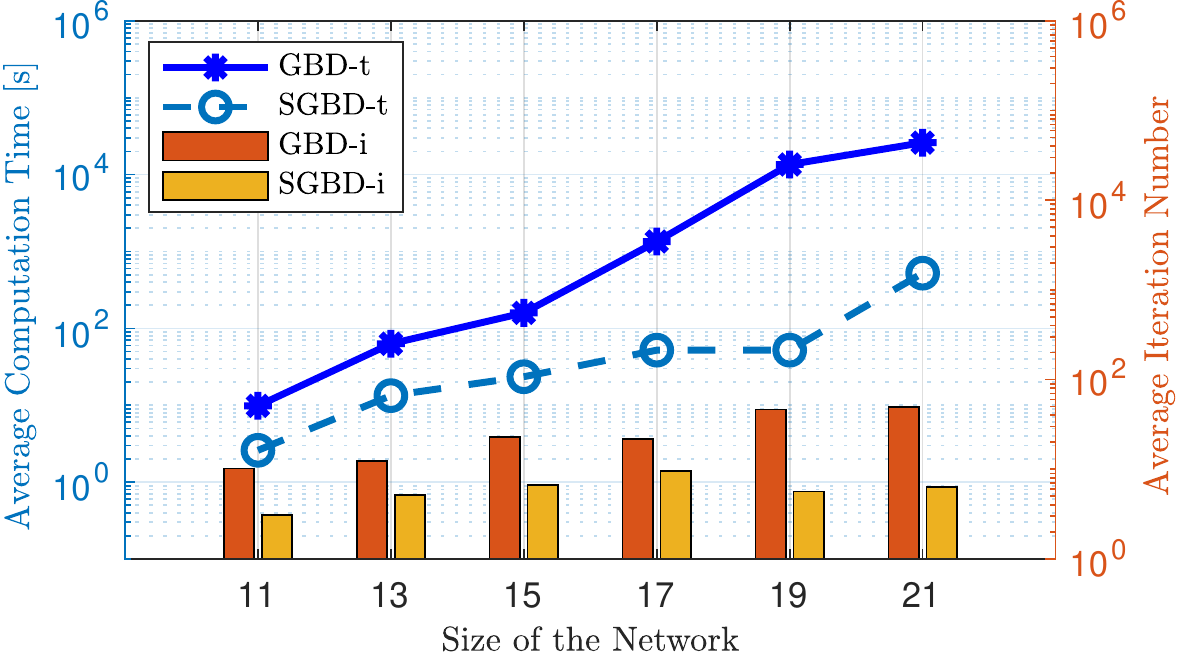}
\caption{Comparison in terms of the average computation time and the average number of iterations. The average computation time is shown on the left Y-axis, whilst the right Y-axis shows the number of iterations. }
\label{Computation_iteration}
\end{figure}

In this subsection, the performance of the GBD and the proposed SGBD  is evaluated in the following two aspects  \text{(i)} the number of iterations, and \text{(ii)} the computation time.  Thanks to the power of the strengthened Benders cuts, the number of iterations of the SGBD (SGBD-i) method will be lower than that of the GBD (GBD-i) method inherently, which has been illustrated in the right Y-axis of Fig.~\ref{Computation_iteration}.  Consequently, compared with the GBD (GBD-t), the SGBD (SGBD-t), which is involved in the computation of the mixed-integer subproblem, still achieves better performances over all instances in terms of the computation time as illustrated in the left Y-axis of Fig.~\ref{Computation_iteration}.



\subsubsection{The scalability of SGBD}
To evaluate the scalability of SGBD, the numerical experiments are conducted with the number of nodes increasing from 11 to 41. Table~\ref{tab:scalability} shows the average run time, as well as the gap between upper bound and lower bound. The $\epsilon$ of Table~\ref{tab:scalability} denotes the gap between upper bound and lower bound, and $*$ from the third row indicates the computation time reaching the time limit of 2 hours. 




\begin{table}[ht]
            \centering
        \caption{The scalability of SGBD}
    \begin{tabular}{ccccccc}
\hline\hline
Instance &  $|\mathcal{E}|=11$   & $|\mathcal{E}|=21$  & $|\mathcal{E}|=31$  & $|\mathcal{E}|=41$      \\
 \hline\hline
Run time~($s$)  &1.53 &   15.12    & 553.54 & $*$	   \\
\rowcolor{Gray}
$\epsilon$   & 0   &  0	& 0\%  & 12\%   \\
 \hline\hline
         \end{tabular}
\label{tab:scalability}
\end{table}



\section{Conclusion}\label{sec:conclusion}
In this paper, we presented a bi-level optimization model to formulate the EVRP with flexible time windows, incentive-aware customers, and a spatio-temporal varying electricity price. Besides, we proposed an accurate reformulation method,  and a strengthened generalized Benders decomposition method combining the complementary advantages of the Benders decomposition method and the Lagrangian dual decomposition approach, to reformulate the original bi-level optimization model and solve the resulting single-level optimization problem, respectively. The simulation results demonstrated that allowing for the time flexibility and jointly optimizing the monetary incentives offered to the customers in exchange of the flexible time window can improve the overall operational cost of the fleet operator by up to  $5\%$, at the same time customers could reduce over $13\%$ of the entire delivery fee. Particularly, a larger upper bound of the time flexibility $\bar{\delta}$ and/or the time-sensitivity $\gamma_2$ specified by the customers surely leaded to the larger delivery fee saving for customers. In addition, even in the most extreme scenarios, we showed that it is always beneficial for the fleet operator to offer monetary incentives, which leverages the flexibility of the resulting delivery time windows.



\section*{Acknowledgments}
The authors thank Dr.\ I.\ New., F. Paparella, and O. Borsboom for proofreading this paper.

\appendix

\section{Proof of Lemma 1}\label{proof:lemma}
\allowdisplaybreaks
 Without loss of generality, we use $\mathcal{X}\coloneqq \{(X_\mathrm{c},Z_\mathrm{d})| \eqref{ssp:a},\eqref{ssp:b},\eqref{ssp:c}, X_\mathrm{c}\in\mathbb{R}^n \}$. When the following condition holds for any $(X_\mathrm{d},\Theta)$,
\begin{equation}\notag
    \Theta \geq \min\limits_{\mathcal{X}\cup Z_\mathrm{d}\in\mathbb{B}^m }    \left\{c_\mathrm{c}^T X_\mathrm{c}, Z_\mathrm{d} = X_\mathrm{d}  \right\},
\end{equation}
the strengthened Benders optimality cut \eqref{cuts:ssp} is valid.
Given any $(X_\mathrm{d},\Theta)$ satisfying the above inequality, we have
\begin{equation}\notag
\begin{aligned}
    \Theta &\geq \min\limits_{\mathcal{X}\cup Z_\mathrm{d}\in\mathbb{B}^m }\quad   c_\mathrm{c}^T X_\mathrm{c}
    \\&\geq \max\limits_{\zeta\in\mathbb{R}^m}\left\{ \zeta^TX_\mathrm{d} +\min\limits_{\mathcal{X}\cup Z_\mathrm{d}\in\mathbb{B}^m }(c_\mathrm{c}^T X_\mathrm{c}-\zeta^T Z_\mathrm{d}  )    \right\}
    \\& = \max\limits_{\zeta\in\mathbb{R}^m}\left\{ \zeta^T (X_\mathrm{d} - X^*_\mathrm{d}) +\min\limits_{\mathcal{X}\cup Z_\mathrm{d}\in\mathbb{B}^m }c_\mathrm{c}^T X_\mathrm{c}-\zeta^T (Z_\mathrm{d} - X^*_\mathrm{d})      \right\}
    \\& = \max\limits_{\zeta\in\mathbb{R}^m}\left\{ \zeta^T (X_\mathrm{d} - X^*_\mathrm{d}) + c_\mathrm{c}^T \bar{X}_\mathrm{c}-\zeta^T (\bar{Z}_\mathrm{d} - X^*_\mathrm{d})      \right\}
    \\& = \max\limits_{\zeta\in\mathbb{R}^m}\left\{ \zeta^T (X_\mathrm{d} - \bar{Z}_\mathrm{d}) + c_\mathrm{c}^T \bar{X}_\mathrm{c}    \right\}
    \\& \geq \zeta^T (X_\mathrm{d} - \bar{Z}_\mathrm{d}) + c_\mathrm{c}^T \bar{X}_\mathrm{c},   
\end{aligned}
\end{equation}
where the second line follows from weak duality of the Lagrangian dual problem and the fourth row follows from the optimality of $(\bar{X}_d, \bar{X}_\mathrm{c})$.  Thus, the strengthened Benders optimality cut is valid.

\section{Proof of Theorem 1}\label{proof:theorem}
Given the solution of Problem~\eqref{BMP} $X_\mathrm{d}^*$, and $\bar{\zeta}\in\mathbb{R}^m$,
\begin{equation}\notag
\begin{aligned}
    \Theta &\geq \max\limits_{\zeta\in\mathbb{R}^m}\left\{ \zeta^T X^*_\mathrm{d} +\min\limits_{\mathcal{X}\cup Z_\mathrm{d}\in\mathbb{B}^m }(c_\mathrm{c}^T X_\mathrm{c}-\zeta^T Z_\mathrm{d}  )    \right\}
    \\& \geq \bar{\zeta}^T X^*_\mathrm{d} +\min\limits_{\mathcal{X}\cup Z_\mathrm{d}\in\mathbb{B}^m }(c_\mathrm{c}^T X_\mathrm{c} - \bar{\zeta}^T Z_\mathrm{d}  )    
    \\& \geq \bar{\zeta}^T X^*_\mathrm{d} +\min\limits_{\mathcal{X}\cup Z_\mathrm{d}\in\mathbb{R}^m  }(c_\mathrm{c}^T X_\mathrm{c} - \bar{\zeta}^T Z_\mathrm{d}  ). 
\end{aligned}
\end{equation}

Note that the second and third inequalities correspond to the strengthened and classic generalized Benders optimality cuts, respectively. These two optimality cuts are parallel as they have the same slope. Compared with the optimality cut obtained from solving the generalized Benders decomposition method, the tightness of the strengthened Benders optimality cut can be quantified by
 \begin{equation}\notag
	\Xi =  \min\limits_{\mathcal{X}\cup Z_\mathrm{d}\in\mathbb{B}^m  }(c_\mathrm{c}^T X_\mathrm{c} - \bar{\zeta}^T Z_\mathrm{d}  )   - \min\limits_{\mathcal{X}\cup Z_\mathrm{d}\in\mathbb{R}^m  }(c_\mathrm{c}^T X_\mathrm{c} - \bar{\zeta}^T Z_\mathrm{d}  ).    
\end{equation} 
 Due to the positive value of $\Xi$, the derived optimality cut is $\Xi$ tighter than the classic generalized optimality cut.


\section{Electric Vehicle Routing Problem without Incentives}\label{app:noincentives}

 \begin{equation}\label{EVRP}\notag
\begin{aligned}
\min\limits_{x^k_{i j},B_{i\tau}, B^\mathrm{s}_{i\tau}, r_i^k, q_j , t_j }  \sum\limits_{i\in\mathcal{R}}\sum\limits_{\tau\in\Lambda}  & \frac{ p_{i\tau} B_{i\tau}\Delta\tau}{g_i} +\sum\limits_{k\in \mathcal{K}}\sum\limits_{i\in\mathcal{V}} \sum\limits_{j\in\mathcal{V}}   (c_{i}   
 + \omega_\mathrm{T} T_{ij}  +  \omega_\mathrm{T}  r^k_{i} g_{i}  ) x^k_{i j}
  \end{aligned}
\end{equation}
\begin{subequations}
    \begin{align}
        \text{s.t.~}  &\sum\limits_{j\in\mathcal{V}}  x^k_{i j}-\sum\limits_{j\in\mathcal{V}} x^k_{j i}=b_{i}, \quad \forall  i \in \mathcal{V}; k\in\mathcal{K},\quad b_{v_1}=1, b_{v_n}=-1, b_{i}=0,
        \\&\sum\limits_{k\in\mathcal{K}}\sum\limits_{j\in\mathcal{V}} x^k_{ij} \leq 1,\quad\forall i\in\mathcal{R} \\
		&t_{j}   \geq T_{ij}+g_ir^k_{i}+t_{i} -M(1-x^k_{ij}),  \forall i\in \mathcal{C}, \quad j \in \mathcal{V}\setminus  v_1,k\in\mathcal{K}\\
		 &t_{j}   \geq T_{ij} + t_{i} -M(1-x^k_{ij}),  \forall i\in \mathcal{R},\quad j \in \mathcal{V}\setminus  v_1,k\in\mathcal{K}	
	\\&t^\mathrm{L}_{j} \leq	t_{j} \leq  t^\mathrm{U}_{j} ,  \forall  j \in \mathcal{R}
	\\&-M (1-x^k_{ij})\leq -E^k_j+E^k_{i}+r_i^k-e_{i j} x^k_{i j}  \leq M (1-x^k_{ij}),\forall i\in\mathcal{C}, j\in\mathcal{V}\setminus  v_1,k\in\mathcal{K}
		\\&-M (1-x^k_{ij})\leq -E^k_j+E^k_{i}-e_{i j} x^k_{i j}  \leq M (1-x^k_{ij}),\forall i\in\mathcal{R}, j\in\mathcal{V}\setminus  v_1,k\in\mathcal{K}
	\\& 0\leq E^k_i\leq E_{k,\mathrm{max}}, \quad  i\in \mathcal{V},k\in\mathcal{K}    
	\\& E^k_{v_1}=E^k_0, \quad\forall k\in\mathcal{K}.\\
	   & c_i=\left\{
		\begin{array}{ll}
			{D_i,} & {\text { if } i\in\mathcal{R}} \\ 
			{c_{v},} & {\text { if } i=v_1}
		\end{array}
		\right.
		  \\& B^\mathrm{s}_{i\tau} \geq B_{i\tau} -B_{i(\tau-1)},\forall i\in\mathcal{C},\tau\in\Lambda \\
  &     \sum\limits_{\tau\in\Lambda} B^\mathrm{s}_{i\tau}\leq 1,\forall i\in\mathcal{C}\\
  &  \sum\limits_{\tau \in \Lambda} B^s_{i \tau} \tau \Delta \tau \leq t_{i}, \quad \forall i \in \mathcal{C} \\
  & t_{i} \leq \sum\limits_{\tau \in \Lambda} B^s_{i \tau}(\tau+1)\Delta \tau, \quad \forall i \in \mathcal{C}\\
  &r^k_i g_i \leq \sum\limits_{\tau\in\Lambda} B_{i\tau} \Delta\tau .
    \end{align}
\end{subequations}

\newpage
\bibliography{references}

\end{document}